\documentclass[titlepage,11pt]{article}
\usepackage{latexsym}
\usepackage{amsfonts}
\usepackage{amsmath}
\newcommand\blackslug{\hbox{\hskip 1pt \vrule width 4pt height 8pt depth 1.5pt
        \hskip 1pt}}
\newcommand\bbox{\hfill \quad \blackslug \bigbreak}

\def\c{\hbox{-}\cdots\hbox{-}}

\def\l{,\ldots,}
\def\ll{,\ldots,}

\title{Induced subgraphs of graphs with large chromatic number.\\
VI. Banana trees}
\author{Alex Scott\\
Mathematical Institute, University of Oxford, Oxford OX2 6GG, UK
\\
\\
Paul Seymour\thanks{Supported by ONR grant N00014-14-1-0084 and NSF
grant DMS-1265563.}\\
Princeton University, Princeton, NJ 08544}

\date{February 21, 2016; revised \today}

\newtheorem{thm}{}[section]

\newcommand{\Proof}{\noindent{\bf Proof.}\ \ }

\begin{document}
\maketitle
\begin{abstract}
We investigate which graphs $H$ have the property that in every graph with bounded clique number and 
sufficiently large chromatic number, some induced subgraph is isomorphic to a subdivision of $H$.
In an earlier paper~\cite{scott}, the first author proved that every tree has this property; and in another earlier paper
with Maria Chudnovsky~\cite{longholes}, we proved that every cycle has this property. Here we give a common generalization. Say a 
``banana'' is the union of a set of paths all with the same ends but otherwise disjoint. We prove that if $H$
is obtained from a tree by replacing each edge by a banana then $H$ has the property mentioned.

\end{abstract}

\section{Introduction}
All graphs in this paper are finite and simple. For some purposes it is convenient to use multigraphs instead of graphs; 
all multigraphs in this paper
are finite and loopless.
If $G$
is a graph, $\chi(G)$ denotes its chromatic number, and $\omega(G)$ denotes its clique number, that is, the cardinality
of the largest clique of $G$.

Let $H$ be a multigraph, and let $J$ be a graph obtained from $H$ by replacing each edge $uv$ by a path (of length at least one)
joining $u,v$, such that these paths are vertex-disjoint except for their ends. 
Then $J$ is a {\em subdivision} of $H$.
We say a graph $G$ is {\em $H$-subdivision-free} if no induced subgraph of $G$ is a subdivision of $H$. 
We could ask: 
\begin{itemize}
\item which multigraphs $H$ have the property that for all $\kappa$ there exists $c$ such that every $H$-subdivision-free graph with clique number
at most $\kappa$ has chromatic number at most $c$?
\item which multigraphs $H$ have the property that for every subdivision $J$ of $H$ and for all $\kappa$ 
there exists $c$ such that every $J$-subdivision-free graph with clique number
at most $\kappa$ has chromatic number at most $c$?
\end{itemize}
The second question, while more complicated, is perhaps better. At least if we confine ourselves to ``controlled'' classes of graphs
(defined later), we know the answer to the second question, while the first remains open. 
The second question could be rephrased as asking for which multigraphs $H$ every graph with bounded clique
number and large chromatic number contains a ``long'' subdivision of $H$, that is, one in which every edge is subdivided at 
least some prescribed number of times.

Let us say a multigraph $H$
is {\em pervasive} in some class of graphs $\mathcal{C}$ if it has the second property above for graphs in the class; that is,
for every subdivision $J$ of $H$ and for all $\kappa \ge 0$ 
there exists $c$
such that every $J$-subdivision-free graph $G\in \mathcal{C}$ with $\omega(G)\le \kappa$ satisfies $\chi(G)\le c$.
(The reader is referred to~\cite{strings} for 
a more detailed introduction to the topic of pervasiveness.)

There are some earlier theorems that can be expressed in this language. First, Scott proved that
\begin{thm}\label{treesubd}
{\rm \cite{scott}} Every tree is pervasive in the class of all graphs.
\end{thm}
Second, we proved with Maria Chudnovsky~\cite{longholes} a conjecture of Gy\'arf\'as~\cite{gyarfas} that for all $\kappa,\ell$, every graph
with clique number at most $\kappa$ and sufficiently large chromatic number has an induced cycle of length at least $\ell$; and that can be reformulated as:
\begin{thm}\label{longholes}
{\rm \cite{longholes}} The multigraph with two vertices and two parallel edges is pervasive in the class of all graphs.
\end{thm}
One of our main theorems, and the goal of the first three-quarters of the paper (up to the end of section 6), 
is the following common generalization:
\begin{thm}\label{bananatree}
Let $H$ be a multigraph obtained from a tree by adding parallel edges. Then $H$ is pervasive in the class of all graphs.
\end{thm}
What is known in the converse direction? Chalopin, Esperet, Li and Ossona de Mendez
proved:
\begin{thm}\label{forestofchand}
{\rm \cite{chandeliers}} Every graph that is pervasive in the class of all graphs
is a forest of chandeliers,
\end{thm}
where
\begin{itemize}
\item  a {\em chandelier} is a graph obtained from a tree by adding a new vertex
called the {\em pivot} adjacent to its leaves (we also count the one- and two-vertex complete graphs as
chandeliers, choosing some vertex as pivot);
\item a {\em tree of chandeliers} is either a chandelier or obtained inductively from a
smaller tree of chandeliers by identifying some vertex with the pivot of a new chandelier; and
\item a {\em forest of chandeliers}
is a graph where every component is a tree of chandeliers.
\end{itemize}
A {\em string graph} is the intersection graph of a set of curves in the plane.
The same paper proved a result stronger than \ref{forestofchand}, namely:
\begin{thm}\label{forestofchand2}
{\rm \cite{chandeliers}} Every graph that is pervasive in the class of string graphs
is a forest of chandeliers.
\end{thm}
With M. Chudnovsky, we proved a converse to this: 
\begin{thm}\label{stringgraphs}
{\rm \cite{strings}} Every forest of chandeliers is pervasive in the class of string graphs. 
\end{thm}
The goal of this paper is to investigate pervasiveness in other classes of graphs. Before we go on, 
we need some definitions.

If $X\subseteq V(G)$, $G[X]$ denotes the subgraph induced on $X$, and we write $\chi(X)$
for $\chi(G[X])$. If $\rho\ge 0$ is an integer, then for $v\in V(G)$, $N^{\rho}_G[v]$
means the set of vertices of $G$ with distance at most
$\rho$ from $v$; and $\chi^{\rho}(G)$ denotes the maximum over all vertices $v$ of $\chi(N^{\rho}_G[v])$, or zero for the null 
graph. Usually we speak of ``$G$-distance'' (to mean distance in $G$) rather than just distance, in case there may be some 
ambiguity.
Let us say an {\em ideal} is a class of graphs closed under taking induced subgraphs; that is, a class 
such that for all graphs $G,H$, if $G$ is an induced subgraph of $H$
and $H\in \mathcal{C}$ then $G\in \mathcal{C}$.
(This is sometimes called a ``hereditary class'', but we needed a shorter name.)
If $\mathcal{C}, \mathcal{C'}$ are ideals and $\mathcal{C'}\subseteq \mathcal{C}$, we say $\mathcal{C}'$ 
is a {\em subideal} of $\mathcal{C}$.

An ideal $\mathcal{C}$ is 
\begin{itemize}
\item {\em colourable} if there exists $k$ such that all members of $\mathcal{C}$ have chromatic 
number at most $k$;
\item {\em $\rho$-bounded} (where $\rho\ge 0$ is some integer) if 
there exists $\tau$ such that $\chi^{\rho}(G)\le \tau$ for all $G\in \mathcal{C}$; 
\item {\em $\rho$-controlled} (for some $\rho$) if every $\rho$-bounded subideal of $\mathcal{C}$ is colourable; and
\item {\em controlled} if it is $\rho$-controlled for some $\rho\ge 0$.
\end{itemize}
Roughly, if a graph in a controlled ideal has large chromatic number, then some ball of bounded radius in the graph also has
large chromatic number. Thus, being $\rho$-bounded and being $\rho$-controlled are almost opposite; $\rho$-bounded means 
$\chi^{\rho}(G)$ is bounded for all $G$, and $\rho$-controlled means $\chi^{\rho}(G)$ 
can be chosen as large as we want by choosing $G$ with $\chi(G)$ large.
Controlled ideals have been an important tool in the study of $\chi$-boundedness -- see~\cite{chiboundedsurvey}.

In fact the following significant extension of \ref{stringgraphs} is known:
\begin{thm}\label{controlled}
{\rm \cite{strings}} Every forest of chandeliers is pervasive in every
controlled ideal.
\end{thm}

We would like to know which multigraphs are pervasive in the ideal of all graphs. Subdividing edges in a graph or multigraph
does not change whether the graph or multigraph is pervasive, so it is enough to decide which graphs are pervasive.
(We could have written this paper just working with graphs, but sometimes multigraphs are more convenient.)
Every such graph  
is a forest of chandeliers, so
let $H$ be a forest of chandeliers; in view of the results of \cite{strings}, what do we still need to prove, to show that
$H$ is pervasive in the ideal of all graphs? Let $J$ be a subdivision of $H$, and let $\mathcal{C}$ be the ideal of all
$J$-subdivision-free graphs; we need to show that the members of $\mathcal{C}$ with bounded clique number also have bounded chromatic number.
Suppose not; so for some $\kappa\ge 0$, there is a noncolourable subideal $\mathcal{D}$ of $\mathcal{C}$ such that all 
graphs $G \in \mathcal{D}$ satisfy $\omega(G)\le \kappa$.
In particular, since $H$ is a forest of chandeliers, \ref{controlled} implies that $\mathcal{D}$
is not $\rho$-controlled, for any $\rho$. Let $\rho\ge 0$. Since
$\mathcal{D}$ is not $\rho$-controlled, there is a noncolourable $\rho$-bounded subideal of $\mathcal{D}$.
Thus, a forest of chandeliers $H$ is not pervasive in the ideal of all graphs if and only if for some subdivision $J$ 
of $H$, some $\kappa\ge 0$ and all $\rho\ge 0$
there is a noncolourable $\rho$-bounded ideal of $J$-subdivision-free graphs all with clique number at most $\kappa$. 

Thus we would like to show the negative: that
for every subdivision $J$
of $H$ and all $\kappa\ge 0$ there exists $\rho\ge 0$ such that
every $\rho$-bounded ideal of $J$-subdivision-free graphs $G$ with $\omega(G)\le \kappa$
is colourable. Let us call this being ``weakly widespread''. 

The argument above shows that
a forest of chandeliers $H$ is not pervasive in the ideal of all graphs if and only if it is not weakly widespread.
We know which graphs are pervasive in controlled ideals, and roughly speaking,
the concept of ``weakly widespread'' is the complementary property; a graph is both pervasive in controlled ideals
and weakly widespread if and only if it is pervasive in the ideal of all graphs, which is what we really want to determine.

Here is a slightly stronger property, eliminating $\kappa$. (The reason for using this
strengthening is that it is somewhat simpler, and it is what in fact we proved whenever we have been able to prove the 
weaker property).
Let us say a multigraph $H$ is {\em widespread} if for every subdivision $J$ 
of $H$ there exists $\rho\ge 0$ such that 
every $\rho$-bounded ideal of $J$-subdivision-free graphs $G$ 
is colourable. 
Equivalently, $H$ is widespread if and only if for every subdivision $J$
of $H$ there exists $\rho\ge 0$ such that
for all $\tau\ge 0$ there exists $c\ge 0$ such that every $J$-subdivision-free graph $G$ with
$\omega(G)\le \kappa$ satisfies $\chi(G)\le c$.

Scott~\cite{scott} conjectured that every multigraph is pervasive in the ideal of all graphs, but this was 
disproved by a beautiful construction in \cite{sevenpoles}. Now we have a different question: which
graphs $H$ are widespread? Originally we expected that the answer would be
``if and only if $H$ is a forest of chandeliers'', but ``if'' remains open and ``only if'' turns out to be false;
in the last quarter of this paper we give some 
widespread graphs that are
not forests of chandeliers. So now our best guess is the following resuscitated version of Scott's conjecture:  

\begin{thm}\label{conj}
{\bf Conjecture:} Every multigraph is widespread.
\end{thm}
We are very far from proving this; we still do not know whether
every forest of chandeliers, or indeed every chandelier, is widespread, and conversely, 
all the multigraphs that we have proved to be 
widespread are subdivisions of outerplanar graphs.

With Chudnovsky, we proved that:
\begin{thm}\label{rbounded}
{\rm \cite{strings}} For all $\rho\ge 2$ and every multigraph $J$, every $\rho$-controlled class of $J$-subdivision-free graphs is $2$-controlled.
\end{thm}

Thus \ref{conj} is equivalent to the following, which is nicer (although we do not use \ref{conj2} in this paper):

\begin{thm}\label{conj2}
{\bf Conjecture:} For all graphs $J$ and for all integers $\tau\ge 0$, there exists $c$ such that if $G$ is a graph with 
chromatic number more than $c$, then either some induced subgraph of $G$ is a subdivision of $J$ or
$\chi^2(G)> \tau$.
\end{thm}

If $e=uv$ is an edge of a multigraph $G$, {\em fattening} $e$ means replacing $e$ 
by some nonempty set of parallel edges all joining $u,v$. We will show:
\begin{thm}\label{fattree}
Let $T$ be a tree, and let $H$ be a multigraph obtained by fattening the edges of $T$.
Then $H$ is widespread.
\end{thm}
Our first main theorem \ref{bananatree} 
follows immediately from \ref{fattree} and \ref{controlled}, as fattening the edges of a tree gives a multigraph whose 
subdivisions are banana trees, and banana trees are trees of chandeliers.

Similar methods can be used to prove that some other classes of multigraphs are widespread.
In sections 7 and 8 we will prove:
\begin{thm}\label{cycle}
Let $H$ be a multigraph obtained from a cycle by fattening all of its edges except one. Then $H$ is    
widespread.
\end{thm}
\begin{thm}\label{triangle}
Let $H$ be a multigraph obtained from a triangle $K_3$ by fattening two of its edges and replacing the third by two parallel edges. 
Then $H$ is widespread.
\end{thm}
These two results are of particular interest because the multigraphs $H$ of \ref{cycle} and \ref{triangle} are in general 
not forests of chandeliers. 

\section{Distant subgraphs with large chromatic number}

With \ref{conj} in mind, let us see what we need. We have a graph $J$ (a subdivision of the initial multigraph $H$), 
and we need to show that
if we choose $\rho$ large enough, then every $\rho$-bounded ideal of $J$-subdivision-free graphs 
is colourable. At this stage we 
prefer not to specify $H,J$, and see how far we can progress in general. So $H,J$ might as well both be 
$K_{\nu}^1$, the graph obtained from a complete graph $K_{\nu}$ by subdividing every edge once; because for any fixed graph $H$,
 if $\nu$ is large enough then
there is an induced subgraph of $K_{\nu}^1$ which is a subdivision of $H$. So, we are given $\nu$, and let us choose $\rho$
very large in terms of $\nu$. 
Now we need to show that every $\rho$-bounded ideal $\mathcal{C}$ of $K_{\nu}^1$-subdivision-free graphs is colourable.
Choose some such $\mathcal{C}$; then since it is $\rho$-bounded, 
there exists $\tau$ such that $\chi^{\rho}(G)\le \tau$ for all $G\in \mathcal{C}$. Altogether then we have three numbers
$\nu, \rho, \tau$, where $\nu$ is given, and $\rho$ is some large function of $\nu$ 
that we can choose, 
and then after selecting $\rho$, the number $\tau$ is given. We need to prove for such a quadruple of numbers, there is a number $c$, 
such that every graph $G$ that is $K_{\nu}^1$-subdivision-free and satisfies $\chi^{\rho}(G)\le \tau$ also satisfies
$\chi(G)\le c$. 
We begin by proving some lemmas about such graphs $G$.
The main result of this section is:

\begin{thm}\label{moredistant}
For all $\nu,k,d, c,\tau\ge 0$ there exists $c'\ge 0$ with the following property.
Let $G$ be a $K_{\nu}^1$-subdivision-free graph, such that $\chi^{2d+7}(G)\le \tau$,
and let $Z\subseteq V(G)$ with $\chi(Z)>c'$. Then there exist subsets
$Z_1\l Z_k\subseteq Z$ such that $\chi(Z_i)>c$ for $i = 1\l k$ and the $G$-distance between every two of $Z_1\l Z_k$ is more than $d$.
\end{thm}

To prove this, we need first:

\begin{thm}\label{getstar}
For all $\kappa\ge 0$ and $d,s\ge 0$ there exists $k\ge s$ with the following property.
Let $G$ be a connected graph with $\omega(G)\le \kappa$. Let $x_1\l x_k\in V(G)$ be distinct, and let $v\in V(G)$ 
such that the $G$-distance between
$v$ and $x_i$ is at most $d$ for $1\le i\le k$. Then there exist $u\in V(G)$ and a subset $S\subseteq \{1\l k\}$ with $|S| = s$,
and for each $i\in S$ a path $Q_i$ between $u$ and $x_i$ of length at most $d$, such that 
\begin{itemize}
\item for $1\le i\le k$, $Q_i$ is a shortest path in $G$ 
between $u,x_i$, and $Q_1\ll Q_k$ all have the same length; 
\item for all distinct $i,j$, $Q_i,Q_j$
are disjoint except for $u$, and there is no edge between $V(Q_i)\setminus \{u\}$ and 
$V(Q_j)\setminus \{u\}$.
\end{itemize}
\end{thm}
\Proof For fixed $\kappa,s$ we proceed by induction on $d$. The result is vacuous for $d=0$ and true for $d=1$, so we assume $d>1$
and that setting $k=k'$ satisfies the result for $d-1$. We may therefore assume that
\\
\\
(1) {\em For each $u\in V(G)$ there are at most $k'-1$ values of $i$ such that $u$ has $G$-distance at most $d-1$ from $x_i$.}
\\
\\
Let $k_1\ge 2$ be such that every graph with at least $k_1$ vertices either has a clique of cardinality $\kappa$ or a stable set of
cardinality $s$, and let $k = 2((d+1)(k'-1)+1)k_1$. We claim that $k$ satisfies the conclusion of the theorem. For given $G,x_1\l x_k,v$ as in the theorem,
let $Q_i$ be a shortest path in $G$ between $v,x_i$ for $1\le i\le k$; thus each $Q_i$ is an induced path. 
Let $D$ be the digraph with vertex set $\{1\l k\}$ in which
for distinct $i,j\in \{1\ll k\}$, there is an edge from $i$ to $j$ if some vertex of $Q_i$ has $G$-distance 
less than $d$ from $x_j$. Because of (1), 
this digraph has outdegree
at most $(d+1)(k'-1)$, and so its underlying graph is $2(d+1)(k'-1)$-degenerate and therefore $(2(d+1)(k'-1)+1)$-colourable; 
and so there exists
$I\subseteq \{1\l k\}$ with
$|I|\ge k/(2(d+1)(k'-1)+1)=k_1$ such that for all distinct $i,j\in I$, every vertex of $Q_i$ has $G$-distance at least $d$ from $x_j$.
It follows that
\begin{itemize}
\item all the paths $Q_i\;(i\in I)$ have length exactly $d$ (to see this, choose $j\in I\setminus \{i\}$; then $v\in V(Q_j)$ 
and so has $G$-distance $d$ from $x_i$);
\item all the paths $Q_i\;(i\in I)$ are pairwise disjoint except for $v$; and
\item for all distinct $i,j\in I$, every edge between $V(Q_i)\setminus \{v\}$ and $V(Q_j)\setminus \{v\}$ joins two neighbours of~$v$.
\end{itemize}
From the choice of $k_1$, since $\omega(G)\le \kappa$, there exists $J\subseteq I$ with $|J|=s$ such that the neighbours
of $v$ in $Q_i\;(i\in J)$ are pairwise nonadjacent; and so the paths $Q_i\;(i\in J)$ satisfy the conclusion of the theorem. 
This proves \ref{getstar}.~\bbox

\begin{thm}\label{router}
For all $\kappa,\nu,d\ge 0$, there exist $k,\ell\ge 0$ with the following property. Let $G$ be a graph with $\omega(G)\le \kappa$, let
$X_1\l X_k$ be nonnull connected subgraphs of $G$, and let $v_1\l v_{\ell}\in V(G)$, such that
\begin{itemize}
\item for all distinct $i,j\in \{1\l k\}$, every vertex in $X_i$ has $G$-distance at least three from every vertex of $X_j$;
\item for all distinct $i,j\in \{1\l \ell\}$, the $G$-distance between $v_i, v_j$ is at least $2d+2$; and
\item for $1\le i\le k$ and $1\le j\le \ell$, the $G$-distance between $X_i$ and $v_j$ is at most $d$.
\end{itemize}
Tnen $G$ is not $K_{\nu}^1$-subdivision-free.
\end{thm}
\Proof 
Let $s=\nu(\nu-1)/2$; choose $k_1$ such that setting $k=k_1$ satisfies \ref{getstar}; let $\ell_1= \nu \binom{k_1}{s}$;
and define $k=(d+1)k_1$ and $\ell=(d+1)^k\ell_1$. We claim that $k,\ell$ satisfy the theorem.
For let $G, X_1\l X_k, v_1\l v_{\ell}$ be as in the theorem.

For $1\le i\le k$ and $1\le j\le \ell$, let the shortest path in $G$ between $X_i$ and $v_j$ have length $d_{ij}$. For
each $j\in \{1\ll \ell\}$, there are only $(d+1)^k$ possibilities for the sequence $(d_{1j}\l d_{kj})$, and so there exists 
$J_1\subseteq \{1\l \ell\}$ with $|J_1|\ge \ell(d+1)^{-k} = \ell_1$ such that for each $i\in \{1\l k\}$, the numbers
$d_{ij}(j\in J_1)$ all have some common value, say $d_i$. Since there are only $d+1$ possibilities for  $d_i$, 
there exists $I_1\subseteq\{1\l k\}$ with $|I_1|=k/(d+1)=k_1$ such that the numbers $d_i\;(i\in I_1)$ all have some common value, say $D$. Thus for each $i\in I_1$
and each $j\in J_1$, the $G$-distance between $X_i$ and $v_j$ is $D$. For each $i\in I_1$ and $j\in J_1$, 
let $P_{ij}$ be some shortest path between $X_i$ and $v_j$, 
and let its end in $X_i$ be $x_{ij}$. 
For each $j\in J_1$, let $G_j$ be the subgraph of $G$ induced on the union of the sets 
$V(P_{ij})\;(i\in I_1)$. For distinct $j,j'\in J_1$, since the $G$-distance between $v_j,v_{j'}$ is at least $2d+2$ and every
vertex of $G_j$ has $G$-distance at most $d$ from $v_j$ and the same for $G_{j'}$, it follows that $G_j,G_{j'}$ are disjoint and there
is no edge joining them. 

Suppose that for some distinct $i,i'\in I_1$ and $j\in J_1$, some vertex $z$ of $P_{ij}$ belongs to or has a neighbour in $X_{i'}$.
Since every path between $x_{ij}$ and $X_{i'}$ has length at least three, it follows that $z$ is not $x_{ij}$ or its neighbour in $P_{ij}$, and so 
there is a path between $v_j$ and $z$ of length at most $D-2$, and hence a path between $v_j$ and $X_{i'}$ of length at most
$D-1$, a contradiction. Thus no vertex of $P_{ij}$ belongs to or has a neighbour in $X_{i'}$.

For each $j\in J_1$, by \ref{getstar} applied to $G_j$ there exist $I^j\subseteq I_1$ with $|I^j|=s$, and a vertex 
$u_j\in V(G_j)$, and induced paths $Q_{ij}$ of $G_j$ between $u_j$ and $x_{ij}$ for each $i\in I^j$, such that 
the paths $Q_{ij}\;(i\in I^j)$ 
are pairwise disjoint except for $u_j$, and
 there are no edges between them not incident with $u_j$. Since $\ell_1= \nu \binom{k_1}{s}$, there exists $J\subseteq J_1$
with $|J|=\nu$ such that the sets $I^j\;(j\in J)$ are all equal, equal to some $I$ say. 
Since $|I|=s = \nu(\nu-1)/2$,
we can number the members of $I$ as $i_{j,j'}$ where $j,j'\in J$ and $j<j'$. For all $j,j'\in J$ with $j<j'$, let $i=i_{j,j'}$;
the subgraph $Q_{ij}\cup X_i\cup Q_{ij'}$ is connected, and so includes an induced path $R_{jj'}$ of $G$ between $u_j, u_{j'}$.
But then the vertices $u_j\;(j\in J)$ and the paths $R_{jj'}$ provide an induced subgraph isomorphic to a subdivision of $K_{\nu}^1$.
This proves \ref{router}.~\bbox

We deduce the following, which is the main step in the proof of \ref{moredistant}:

\begin{thm}\label{distant}
For all $\nu,d, c,\tau\ge 0$ there exists $c'\ge 0$ with the following property.
Let $G$ be a $K_{\nu}^1$-subdivision-free graph, such that $\chi^{2d+7}(G)\le \tau$, 
and let $Z\subseteq V(G)$ with $\chi(Z)>c'$. Then there exist subsets 
$Z_1,Z_2\subseteq Z$ such that $\chi(Z_i)>c$ for $i = 1,2$ and the $G$-distance between $Z_1,Z_2$ is more than $d$.
\end{thm}
\Proof Let $k,\ell\ge 1$ satisfy \ref{router} with $d$ replaced by $d+3$, and $\kappa=\tau$. Let $c_k=\ell\tau$, 
and for $i = k-1\l 0$ define $c_i=2c_{i+1}+2c$. Let $c'=c_0$. We claim that $c'$ satisfies the conclusion of the theorem.
For let $G, Z$ be as in the theorem. Then $\omega(G)\le \tau$, since $\omega(G)\le \chi^1(G)\le \chi^{2d+7}(G)\le \tau$.
Choose $k'\le k$ maximum such that there exist connected subgraphs $X_1\l X_{k'}$ of $G[Z]$
and a subset $A\subseteq Z$ with the following properties:
\begin{itemize}
\item for $1\le i<j\le k'$, the $G$-distance between $X_i, X_j$ is at least three;
\item for $1\le i\le k'$, there exists $d_i$ with $3\le d_i\le d+3$ such that 
every vertex in $A$ has $G$-distance exactly $d_i$ from $X_i$; and
\item $\chi(A)>c_{k'}$.
\end{itemize}
(This is possible since setting $k'=0$ and $A=Z$ satisfies the bulletted statements.) By replacing $A$ by
the vertex set of a connected component of $G[A]$ with maximum chromatic number, we may assume 
in addition that $G[A]$
is connected.
Since $c_{k'}> (\ell-1)\tau$ and $\chi^{2d+7}(G)\le \tau$, there exist vertices $v_1\l v_{\ell}\in A$, 
pairwise with $G$-distance at least $2d+8$. Consequently, $k'<k$ by \ref{router}.
Let $z_0\in A$. For $i\ge 0$ let $L_i$
be the set of vertices in $A$ with $G[A]$-distance from $z_0$ equal to $i$. For $i\ge 0$ let $M_i = L_0\cup\cdots\cup L_i$.
Thus each $M_i$ induces a connected graph.
For $r\ge 0$ let $M_i^r$ denote the set of vertices in $A$ with $G$-distance from $M_i$ at most $r$.
It follows that $M_i^r \supseteq M_{i+r}$.
(We emphasize that $G$-distance and $G[A]$-distance may be different.)
For sufficiently large $i$, $M_i = A$; and so there exists $i$ such that $\chi(M_i^{d+3})>2c_{k'+1}+c$. 
Choose $i$ minimum with this property.

Suppose that $\chi(M_i^2)\le c$. Then $\chi(M_i^{d+3}\setminus M_i^2) > (2c_{k'+1}+c)-c$, and every vertex in $M_i^{d+3}\setminus M_i^2$
has $G$-distance from $M_i$ at least three and at most $d+3$. For $3\le j\le d+3$ let $B_j$ be the set of vertices in
$M_i^{d+3}\setminus M_i^2$ with $G$-distance exactly $j$ from $M_i$. It follows that 
$\chi(B_j)>c_{k'+1}$ for some $j\in \{3\l d+3\}$. Let $X_{k'+1} = M_i$ and $d_{k'+1} = j$; then since $\chi(B_j)> c_{k'+1}$,
this contradicts the maximality of $k'$. This proves that $\chi(M_i^2)>c$.

Now $\chi(M_i^{d+3})>2c_{k'+1}+c\ge \tau$, and so $i>0$ since $\chi^{2d+7}(G)\le \tau$.
From the minimality of $i$ it follows that $\chi(M_{i-1}^{d+3})\le 2c_{k'+1}+c$. Since $\chi(A)>c_{k'}$, it follows that 
$\chi(A\setminus M_{i-1}^{d+3})>c_{k'}-(2c_{k'+1}+c)= c$.
But $M_i^2\subseteq M_{i-1}^3$, so the $G$-distance between $M_i^2$ and $A\setminus M_{i-1}^{d+3}$ is at least $d+1$.
Since both the sets are subsets of $A$ and hence of $Z$, and both sets have chromatic number more than $c$, 
this proves \ref{distant}.~\bbox

Now we can prove the main result of this section \ref{moredistant}, which we restate:
\begin{thm}\label{moredistant2}
For all $\nu,k,d, c,\tau\ge 0$ there exists $c'\ge 0$ with the following property.
Let $G$ be a $K_{\nu}^1$-subdivision-free graph, such that $\chi^{2d+7}(G)\le \tau$,
and let $Z\subseteq V(G)$ with $\chi(Z)>c'$. Then there exist subsets
$Z_1\l Z_k\subseteq Z$ such that $\chi(Z_i)>c$ for $i = 1\l k$ and the $G$-distance between every two of $Z_1\l Z_k$ is more than $d$.
\end{thm}
\Proof We proceed by induction on $k$. Let $c''$ satisfy the theorem with $k$ replaced by $k-1$; and let $c'$ satisfy
\ref{distant} with $c$ replaced by $c''$. We claim that $c'$ satisfies the conclusion of the theorem. For let $G,Z$ be as in the theorem.
By \ref{distant} exist subsets
$Z_1,Z_2\subseteq Z$ such that $\chi(Z_i)>c''$ for $i = 1,2$ and the $G$-distance between $Z_1,Z_2$ is more than $d$.
By the inductive hypothesis applied with $Z$ replaced by $Z_2$, there are $k-1$ subsets $Y_1\l Y_{k-1}$ of $Z_2$, each with
chromatic number at least $c$ and pairwise with $G$-distance at least $d+1$. But then $Z_1,Y_1\l Y_{k-1}$ satisfy the theorem. This proves \ref{moredistant2}~\bbox

\section{Pineapple trees}

If $X,Y\subseteq V(G)$, we say that $Y$ {\em covers} $X$ if $X\cap Y=\emptyset$ and every vertex in $X$ has a neighbour in $Y$.
If in addition $G[Y]$ is connected we call the pair $(X,Y)$ a {\em pineapple} in $G$. It is a {\em levelled} pineapple
if there exists $z_0\in Y$ such that for some $k$, every vertex in $Y$ is joined to $z_0$ by a path of $G[Y]$ of length less than $k$,
and there is no path in $G[X\cup Y]$ of length less than $k$ from $z_0$ to $X$.

Now let $T$ be a tree, with a vertex $r$ called its {\em root}. We call $(T,r)$ a {\em rooted tree}.
For $u,v\in V(T)$, we say $v$ is an {\em ancestor} of $u$ and $u$ is a {\em descendant} of $v$
if $v$ belongs to the path of $T$ between $u,r$. We define {\em parent} and {\em child} in the natural way. 
We say $u,v\in V(T)$ are {\em incomparable} if neither is a descendant of the other. Let $L(T)$
be the set of vertices of $T$ with no children (thus, $L(T)$ is the set of leaves of $T$ different from $r$, 
except when $V(T)=\{r\}$). Now let $G$ be a graph.
For each vertex $v\in V(T)$ let $C_v\subseteq V(G)$, and for each vertex $v\in V(T)\setminus L(T)$
let $(X_v,Y_v)$ be a levelled pineapple in $G$ with $X_v\cup Y_v=C_v$, and with the following properties:
\begin{itemize}
\item all the sets $C_v\;(v\in V(T))$ are nonempty and pairwise disjoint;
\item for all incomparable $u,v\in V(T)$ there is no edge between $C_u,C_v$;
\item if $u, v\in V(T)$ are distinct, and $u$ is a descendant of $v$, then there are no edges between
$C_u$ and $Y_v$, and if also $u\in L(T)$ then $X_v$ covers $C_u$.
\end{itemize}
(Note that we only demand that $X_v$ covers $C_u$ when $u$ is a 
leaf. It is undetermined
whether there are edges between $C_u$ and $X_v$ when $u\in V(T)\setminus L(T)$ is a descendant of $v$; this will be resolved later.)
We call the system 
$$(T,r,((X_v,Y_v):v\in V(T)\setminus L(T)), (C_v:v\in V(T)))$$ 
a {\em pineapple tree} in $G$ and $(T,r)$ is its {\em shape}. 
Let us call the union of all the sets $C_v(v\in V(T))$
the {\em vertex set} of the pineapple tree.
In this section we prove:

\begin{thm}\label{gettree}
For all $\nu,c,d, \tau\ge 0$, and every rooted tree $(T,r)$, there exists $c'$ with the following property.
Let $G$ be a $K_{\nu}^1$-subdivision-free graph, such that $\chi^{2d+7}(G)\le \tau$, and let $Z\subseteq V(G)$
with $\chi(Z)>c'$.
Then there is a pineapple tree
$$(T,r,((X_v,Y_v):v\in V(T)\setminus L(T)), (C_v:v\in V(T)))$$
in $G$, with vertex set a subset of $Z$, such that
\begin{itemize}
\item $\chi(C_v)>c$ for each $v\in L(T)$; and 
\item for all incomparable $u,v\in V(T)$, 
the $G$-distance between $C_u,C_v$
is at least $d+1$.
\end{itemize}
\end{thm}
\Proof First, a remark: the reader might wonder why we do not just delete all vertices in $V(G)\setminus Z$. 
The reason is, the final condition involves $G$-distance and not just $G[Z]$-distance. 

We proceed by induction on $|V(T)|$. If 
$V(T)=\{r\}$, then $r\in L(T)$, and we define $C_r=V(G)$ and the theorem holds. Thus we may assume that $r\notin L(T)$, and the result
holds for smaller trees. 
Let $r_1\l r_k$ be the children of $r$, and for $1\le i\le k$ let $T_i$ be the component of $T\setminus r$ containing $r_i$.
Inductively for $i = 1\l k$, there exists $c_i$ satisfying the theorem with $T,r,c'$ replaced by $T_i, r_i, c_i$. Let $c'' $
be the maximum of $c_1\l c_k$. Choose $c_0\ge \tau$ such that \ref{moredistant} holds with $c,c'$ replaced by $c'',c_0$.
We claim that setting $c'=2c_0$ satisfies the conclusion of the theorem. For let $G,Z$ be as in the theorem.
Choose a component $A$ of $G[Z]$ with $\chi(A) = \chi(Z)$, and choose $z_0\in A$. For $j\ge 0$
let $L_j$ be the set of vertices in $A$ with $G[Z]$-distance $j$ from $z_0$, and choose $j$ such that 
$\chi(L_j)\ge \chi(A)/2>c_0$.
It follows that $j>1$, since $\chi^{2d+7}(G)\le \tau\le c_0$; let $X_r = L_{j-1}$ and $Y_r = L_0\cup\cdots\cup L_{j-2}$.
Then $(X_r,Y_r)$ is a levelled pineapple, and $X_r$ covers $L_j$, and there are no edges between $Y_r, L_j$.

From \ref{moredistant} there exist $Z_1\l Z_k\subseteq L_j$, each with chromatic number
more than $c''$, and pairwise at $G$-distance more than $d$. From the choice of $c_i$, for each $i$ there is a pineapple tree
$$(T_i,r_i,((X_v,Y_v):v\in V(T_i)\setminus L(T_i)), (C_v:v\in V(T_i)))$$
in $G$, with vertex set a subset of $Z_i$, such that $\chi(C_v)>c''\ge c_i$ for each $v\in L(T_i)$, 
and for all distinct $u,v\in L(T_i)$, the $G$-distance between $C_u,C_v$
is at least $d+1$. But then 
$$(T,r,((X_v,Y_v):v\in V(T)\setminus L(T)), (C_v:v\in V(T)))$$
is the required pineapple tree. This proves \ref{gettree}.~\bbox

We remark that \ref{gettree} proved more than we asked for; we found a pineapple tree such that $X_v$ covers $C_u$
for all $u,v$ such that $u$ is a descendant of $v$, not just when $u$ is a leaf. This extra property will not be of use to us.

\section{Pruning a pineapple tree}

In the previous section we proved that we can assume our graph contains a pineapple tree with shape whatever we like;
but we
need to extract a nice big banana tree from this somehow. To do so, we will tidy up the existence of
edges between $X_u\cup Y_u$ and $X_v$ when $u\in V(T)\setminus L(T)$ is a descendant of $v$. 
What we care about is, let $r,v,w,u$ be in some path of $T$ in this order, where $u\in L(T)$. 
There are vertices in $X_v$ that have neighbours in
$C_u$; do all such vertices have neighbours in $C_w$? Do none of them have neighbours in $C_w$? 
It turns out that both these extreme cases are good for finding banana trees
(assuming the same thing happens for all such $u,v,w$, and the shape of the pineapple tree 
is rich enough and the chromatic numbers of the sets $C_v$ are big enough). 
So we 
need a Ramsey argument to produce a pineapple tree with one of these two properties, and that is the content of this 
section and the next. We do it in two stages: first we arrange that for each triple $v,w,u$, one of the two things 
happens (this is called ``pruning''), and then we will arrange that the same thing happens for all triples $v,w,u$.

Let $$(T,r,((X_v,Y_v):v\in V(T)\setminus L(T)), (C_v:v\in V(T)))$$ be a pineapple tree in $G$, and let $d\ge 2$, such that
for all incomparable $u,v\in V(T)$, the $G$-distance between $C_u,C_v$
is at least $d+1$. For each $u\in L(T)$ and each ancestor 
$v\in V(T)\setminus L(T)$ of $u$, let 
$X_v^u$ be the set of vertices
in $X_v$ with a neighbour in $C_u$.
Here are some observations
about these subsets:

\begin{itemize}
\item For each $v\in V(T)\setminus L(T)$, if $u,u'\in L(T)$ are distinct descendants of $v$, then $X_v^u\cap X_v^{u'}=\emptyset$;
for $d\ge 2$, so the $G$-distance between $C_u, C_{u'}$ is at least three, and so no vertex in $X_v$ has neighbours in both sets.
\item We may assume that for each $v\in V(T)\setminus L(T)$, every vertex in $X_v$ belongs to $X_v^u$ for some descendant
$u\in L(T)$ of $v$; for any other vertices in $X_v$ may be removed from $X_v$ without violating the definition of a pineapple tree.
\item For all distinct $u,v,v'\in V(T)$, if $u\in L(T)$, and $v$ is an ancestor of $u$, and $v'$ is incomparable with $u$, 
then there are no edges 
between $X_v^u$ and $C_{v'}$; because the $G$-distance between $C_u,C_{v'}$ is at least three, and every vertex in $X_v^u$ has a neighbour
in $C_u$, and so has no neighbour in $C_{v'}$.
\end{itemize}

We say an {\em aligned triple} is a triple $(u,v,w)$ such that $u\in L(T)$, $w$ is an ancestor of $u$, $v$ is an ancestor of $w$, 
and $u,v,w$ are all different. We say that an aligned triple $(u,v,w)$ is {\em pruned} if either
every vertex in $X_v^u$ has a neighbour in $Y_w$, or none does; and the pineapple tree is {\em pruned} if every aligned triple
is pruned.

\begin{thm}\label{pruned}
For all $\nu,c,d, \tau\ge 0$, and every rooted tree $(T,r)$, there exists $c'$ with the following property.
Let $G$ be a $K_{\nu}^1$-subdivision-free graph, such that $\chi^{2d+7}(G)\le \tau$, and let $Z\subseteq V(G)$ 
with $\chi(Z)>c'$.
Then there is a pruned pineapple tree
$$(T,r,((X_v,Y_v):v\in V(T)\setminus L(T)), (C_v:v\in V(T)))$$
in $G$, with vertex set a subset of $Z$, such that
\begin{itemize}
\item $\chi(C_v)>c$ for each $v\in L(T)$; and
\item for all incomparable $u,v\in V(T)$,
the $G$-distance between $C_u,C_v$
is at least $d+1$.
\end{itemize}
\end{thm}
\Proof Let $c''=2^{h^2}c$ where $h$ is the length of the  
longest path in $T$ with one end $r$. Let $c'$ satisfy \ref{gettree} with $c$ replaced by $c''$.
We claim that $c'$ satisfies the conclusion of the theorem. For let $G,Z$ be as in the theorem; then by \ref{gettree} there is a pineapple tree
as in \ref{gettree}, in the usual notation. 
Thus we may choose a pineapple
tree 
$$(T,r,((X_v,Y_v):v\in V(T)\setminus L(T)), (C_v:v\in V(T)))$$
in $G$, satisfying the following conditions:
\begin{itemize}
\item its vertex set is a subset of $Z$;
\item for all incomparable $u,v\in L(T)$, the $G$-distance between $C_u,C_v$
is at least $d+1$;
\item for each $u\in L(T)$, let $n_u$ be the number of pairs $(v,w)$ such that $(u,v,w)$ is a pruned aligned triple;
then $\chi(C_u)>c''2^{-n_u}$ .
\end{itemize}
(Indeed, we can choose such a pineapple tree with  $\chi(C_u)>c''$ for each $i$.)
Choose this tree such that in addition
the sum of the numbers $n_u\;(u\in L(T))$ is maximum. We claim this tree is pruned. For if not, choose an aligned triple $(u,v,w)$
that is not pruned.
Let $A$ be the set of vertices in $X_u^v$ with a neighbour in $Y_w$, and $B=X_v^u\setminus A$.
Every vertex in $C_u$ has a neighbour in one of $A,B$, and so we can choose one of $A,B$, say $W_v^u$, such that the set of vertices in
$C_u$ with a neighbour in $W_v^u$, say $C_u'$, has chromatic number at least $\chi(C_u)/2$ and hence more than $c''2^{-n_u-1}$.
But then replacing $X_v^u$ by $W_v^u$ and $C_u$ by $C_u'$ gives a new pineapple with the sum of the numbers $n_u\;(u\in L(T))$
larger, which is impossible. This proves that the pineapple tree is pruned. 

For each $u\in L(T)$, $n_u\le h(h-1)/2\le h^2$ and so $\chi(C_u)>c''2^{-n_u}\ge c$, and so this pineapple tree satisfies 
the conclusion of 
the theorem.~\bbox

\section{A Ramsey theorem for trees}

Let $h\ge 0$ and $t\ge 1$, and let
$(T,r)$ be a rooted tree in which every path from $r$ to a member of $L(T)$ has length $h$, and every vertex in 
$V(T)\setminus L(T)$ has $t$ children (and hence degree $t+1$, except for $r$). We call $(T,r)$ a
{\em uniform $t$-ary tree of height $h$}. We need the following.

\begin{thm}\label{treeramsey} Let $q,h\ge 0$ and $t\ge 1$.
Let $(T',r)$ be a uniform $(qt)$-ary tree of height $h$, and let $\phi$ be a map from $L(T')$ 
to the set $\{1\l q\}$. Then there is a subtree $T$ of $T'$ containing $r$, such that $(T,r)$ is a 
uniform $t$-ary tree of height $h$, and such that for some $x\in \{1\l q\}$, $\phi(u)=x$ for all $u\in L(T)$.
\end{thm}
\Proof We proceed by induction on $h$. For $h=0$ the result is true, so we assume that $h>0$ and the result holds for $h-1$.
Let $r_1\l r_{qt}$ be the children of $r$ in $T'$, and for $1\le i\le qt$ let $T_i$ be the component of $T\setminus r$
containing $r_i$. For $1\le i\le qt$, from the inductive hypothesis there is a 
subtree $T_i'$ of $T_i$ containing $r_i$, such that $(T_i',r_i)$ is a
uniform $t$-ary tree of height $h-1$, and such that for some $x_i\in \{1\l q\}$, $\phi(u)=x_i$ for all $u\in L(T_i')$.
Choose $x\in \{1\l q\}$ such that $x_i=x$ for at least $t$ values of $i$; then the union of $t$ of the corresponding trees $T_i$,
together with $r$, gives the desired tree $T'$.~\bbox

Let
$$(T,r,((X_v,Y_v):v\in V(T)\setminus L(T)), (C_v:v\in V(T)))$$
be a pineapple tree. It is {\em barren} if for every aligned triple $(u,v,w)$, no member of $X_v^u$
has a neighbour in $Y_w$; and it is {\em fruitful} if for every aligned triple $(u,v,w)$, every member of $X_v^u$
has a neighbour in $Y_w$. (In both cases, there may or may not be edges between $X_v^u$ and $X_w$.)
We need a further strengthening of \ref{gettree}.
\begin{thm}\label{platonic}
For all $\nu,c,d, \tau\ge 0$, and every rooted tree $(T,r)$, there exists $c'$ with the following property.
Let $G$ be a $K_{\nu}^1$-subdivision-free graph such that $\chi^{2d+7}(G)\le \tau$, and let $Z\subseteq V(G)$
with $\chi(Z)>c'$.
Then there is a pineapple tree
$$(T,r,((X_v,Y_v):v\in V(T)\setminus L(T)), (C_v:v\in V(T)))$$
in $G$, with vertex set a subset of $Z$, such that
\begin{itemize}
\item $\chi(C_v)>c$ for each $v\in L(T)$; and
\item for all incomparable $u,v\in V(T)$,
the $G$-distance between $C_u,C_v$
is at least $d+1$.
\end{itemize}
which is either barren or fruitful.
\end{thm}
\Proof
Choose $t\ge 1$ and $h\ge 0$ such that every vertex of $T$ has at most $t$ children 
and every path of $T$ with one end $r$ has length at most $h$. Let $q= 2^{2^{2h}}$. Let $(T',r')$ be a uniform 
$(qt)$-ary tree of height $2^h$, and choose $c'$ such that \ref{pruned} is satisfied with $(T,r)$ replaced by $(T'r')$.
We claim that $c'$ satisfies the conclusion of the theorem. For let $G,Z$ be as in the theorem. By \ref{pruned} 
there is a pruned pineapple tree 
$$(T',r',((X_v,Y_v):v\in V(T')\setminus L(T')), (C_v:v\in V(T')))$$ 
as in \ref{pruned}.
For each $u\in L(T')$, let $q_u$ be the function with domain the set of all ordered pairs $(i,j)$ with $0\le i<j< 2^h$, 
defined as follows. For each such pair $(i,j)$, let $v,w\in V(T')$ be the ancestors of $u$ with $G$-distance $i$ and $j$ from $r$
respectively; let $q_u(i,j) = 0$ if no member of $X_v^u$
has a neighbour in $Y_w$, and $q_u(i,j) = 1$ if every member of $X_v^u$
has a neighbour in $Y_w$. (Since the pineapple tree is pruned and all the sets $X_v^u$ are nonempty, this is well-defined.)
Thus each $q_u$ is a map into a domain with at most $q$ elements, and so by \ref{treeramsey} there 
is a subtree $T''$ of $T'$ containing $r'$, such that $(T'',r')$ is a
uniform $t$-ary tree of height $2^h$, and such that all the functions $q_u\;(u\in L(T''))$ are equal.
Let the common value of all the $q_u\;(u\in L(T''))$ be a function $f$. Let $H$ be the graph with vertex set $\{0\l 2^h-1\}$ in which for
$0\le i<j<2^h$, $i,j$ are adjacent if $f(i,j)=1$. By Ramsey's theorem applied to $H$, there exists $I\subseteq \{0\l 2^h-1\}$ with $|I|=h$
such that all the values $f(i,j)\;(i,j\in I, i<j)$ are equal. Let $I=\{i_0\l i_{h-1}\}$ where $0\le i_0<\cdots<i_{h-1}<2^h$.
Choose $s\in V(T'')$ with $T''$-distance $i_0$ from $r'$. Let $N$ be the set of descendants of $s$ in $T''$ 
whose $T''$-distance from $r'$ belongs the set $I\cup \{2^h\}$. 
Let $S$ be the tree with vertex set $N$ in which $u,v$ are adjacent if one is a descendant in $T''$ of the other and no third vertex
of $N$ belongs to the path of $T''$ between them. Then $(S,s)$ is a rooted tree in which every path from $s$ to $L(S)$ has length $h$
and every vertex in $V(S)\setminus L(S)$ has $t$ children. Consequently 
$$(S,s,((X_v,Y_v):v\in V(S)\setminus L(S)), (C_v:v\in V(S)))$$
is a pineapple tree, and it is either barren or fruitful, and since $(S,s)$ has a rooted subtree isomorphic to $(T,r)$, the result follows.
This proves \ref{platonic}.~\bbox

We can eliminate barren pineapple trees, but to do so we need a lemma.
Say an {\em infusion} of a graph $H$ in a graph $G$ is a map $\phi$, such that
\begin{itemize}
\item $\phi$ maps $V(H)$ injectively into $V(G)$,
and
\item $\phi$ maps each edge $e=uv$ of $H$ to an induced path $\phi(e)$ of $G$ between $\phi(u),\phi(v)$, of length at least two;
\item for all distinct $e,f\in E(H)$, the paths $\phi(e),\phi(f)$ are vertex-disjoint except possibly for a common end;
\item for all distinct $e,f\in E(H)$ with no common end, there is no edge of $G$ between $V(\phi(e)),V(\phi(f))$;
\item for all distinct $e,f\in E(H)$ with a common end $v$ say, there is at most one edge between $V(\phi(e)\setminus \phi(v))$
and $V(\phi(f)\setminus \phi(v))$ and such an edge joins the two neighbours of $\phi(v)$.
\end{itemize}

\begin{thm}\label{immersion}
There exists $n$ such that if $G$
is $K_{\nu}^1$-subdivision-free and $\omega(G)\le \kappa$, there is no infusion of $K_n$ in $G$.
\end{thm}
\Proof
The result is an easy application of Ramsey's theorem, and so we merely sketch the proof.
Choose $n$ very large in terms of $\kappa, \nu$, and suppose $G, \phi$ is a counterexample. Let $\mu=\nu^2$.
We may assume that
every vertex of $G$ belongs to one of the paths $\phi(ij)$. Consequently, for each path $\phi(ij)$, all its vertices have 
degree two in $G$, except for the first two and the last two. If $e$ is an edge of $\phi(ij)$ and neither of its ends is one of the first two or last two vertices of $\phi(ij)$, 
then both its ends have degree two in $G$ and we could contract $e$ and make a smaller counterexample; 
and so we may assume that each path $\phi(ij)$ has at most four edges.
Thus we have:
\begin{itemize}
\item the vertices $\phi(v)\;(v\in V(K_n))$ are pairwise nonadjacent (since the paths $\phi(ij)$ all have length at least two
and are induced);
\item each path $\phi(ij)$ has length at most four; and
\item $G$ is $K_{\mu,\mu}^1$-subdivision-free, where $K_{\mu,\mu}^1$ is obtained from the complete bipartite graph $K_{\mu,\mu}$
by subdividing every edge once; because $K_{\mu,\mu}^1$ contains a subdivision of $K_{\nu}^1$ as an induced subgraph.
\end{itemize}
But then the result is a consequence of theorem 3.2 of~\cite{strings}.~\bbox

\begin{thm}\label{easykind}
For all $\nu\ge 0$, there exists $n\ge 1$ with the following property. Let $\tau\ge 0$ and let
$(T,r)$ be a rooted tree where $T$ is a path of length 
$n^2$ with ends $r,u$. Let $G$ be a $K_{\nu}^1$-subdivision-free graph, such that $\chi^{2}(G)\le \tau$.
Then there is no barren pineapple tree in $G$ with shape $(T,r)$, such that in the usual notation $\chi(C_u)>n\tau$.
\end{thm}
\Proof Let $\kappa=\tau$, and choose $n$ as in \ref{immersion}.
We claim that $n$ satisfies the conclusion of the theorem. Let $T,r,G$ be as in the theorem;
since $\omega(G)\le \chi^2(G)\le tau$, it follows that $\omega(G)\le \kappa$.
Suppose that there is a barren pineapple tree with shape $(T,r)$ as described in the theorem. Then $|L(T)|=1$; let $u\in L(T)$. 
Let the vertices of $T$ be $t_0\c t_{n^2}$ in order, where $r = t_0$. (Thus $u= t_{n^2}$.)
For $0\le k<n^2$ let us write $X_k$ for $X_{t_k}$
and $Y_k$ for $Y_{t_k}$ for convenience. Since
$\chi(C_u)>n\tau$ and $\chi^{2}(G)\le \tau$, there exist $n$ vertices $v_1\l v_n$ in $C_u$, pairwise with $G$-distance at least three.
For $1\le i\le n$ and $0\le k< n^2$, let 
$x^i_k$ be a neighbour of $v_i$ in $X_{k}$. Let $H$ be a graph with vertex set $\{v_1\l v_n\}$ in which all pairs of vertices 
are adjacent. Number the edges of $H$ as $e_0\l e_{m-1}$ where $m = n(n-1)/2$. Let $0\le k< m$ and let $e_k$
have ends $v_i,v_j$ say where $i<j$. Let $P_k$ be an induced path of $G$ between $v_i,v_j$, consisting of the edges 
$v_ix^i_k$, $v_jx^j_k$ and an induced path joining $x^i_k, x^j_k$ with interior in $Y_k$ (this exists since $Y_k$ covers $X_k$
and $G[Y_k]$ is connected). Then the paths $P_0\l P_{m-1}$ are pairwise vertex-disjoint except possibly for a common end.
Suppose that there is an edge $e$ of $G$ joining $P_k,P_{k'}$ say, where $k\ne k'$, and $e$ is not incident with a common end
of $P_k, P_{k'}$. Suppose first that some end $v$ of $e$ is an end of one of $P_k, P_{k'}$, say of $P_k$;
then $v=v_i$ for some $i$. Since $v_1\l v_n$ pairwise have $G$-distance at least three, the other end of $e$ is not
an end of $P_{k'}$, and so belongs to $X_{k'}\cup Y_{k'}$. It cannot belong to $Y_{k'}$ since there are no edges between 
$C_u$ and $Y_{k'}$, from the definition of a pineapple tree. Consequently it belongs to $X_{k'}$, and so is adjacent to an end of $P_k'$.
This end of $P_{k'}$ must be $v_i$, since $v_1\l v_n$ pairwise have $G$-distance at least three, and so $v_i$ is a common end of
$P_k, P_{k'}$, contrary to the definition of $e$. Thus neither end of $e$ belongs to $C_u$. Hence one end is in $X_k\cup Y_k$,
and the other in $X_{k'}\cup Y_{k'}$. From the symmetry we may assume that $k'>k$, and so there are no edges between $Y_k$
and $X_{k'}\cup Y_{k'}$ from the definition of a pineapple tree. Hence one end of $e$ is in $X_k$. The other
end of $e$ is not in $Y_{k'}$ since $(u,t_k,t_{k'})$ is an aligned triple, and by hypothesis no vertex in $X_{k}^u$
has a neighbour in $Y_{k'}$. Hence the other end of $e$ is in $X_{k'}$. It follows that both ends of $e$ have a neighbour in 
$\{v_1\l v_k\}$, and so they have a common neighbour since $v_1\l v_n$ pairwise have $G$-distance at least three, 
and this common neighbour
is a common end of both $P_k,P_{k'}$. Consequently the vertices $v_1\l v_n$ and paths $P_0\l P_{m-1}$ define an infusion of $K_n$
in $G$, which is impossible. This proves \ref{easykind}.~\bbox

We deduce:
\begin{thm}\label{summary}
For all $\nu,c,d, \tau\ge 0$, and every rooted tree $(T,r)$, there exists $c'$ with the following property.
Let $G$ be a $K_{\nu}^1$-subdivision-free graph such that $\chi^{2d+7}(G)\le \tau$, and let $Z\subseteq V(G)$
with $\chi(Z)>c'$.
Then there is a fruitful pineapple tree
$$(T,r,((X_v,Y_v):v\in V(T)\setminus L(T)), (C_v:v\in V(T)))$$
in $G$, with vertex set a subset of $Z$, such that 
\begin{itemize}
\item $\chi(C_v)>c$ for each $v\in L(T)$; and 
\item for all incomparable $u,v\in L(T)$, the $G$-distance between $C_u,C_v$
is at least $d+1$.
\end{itemize}
\end{thm}
\Proof By adding a path to $T$ if necessary, we may assume that there is a rooted subtree of $(T,r)$ which is a path
of length $n$ as in \ref{easykind}. But then the result follows from \ref{platonic} and \ref{easykind}.~\bbox

\section{Banana trees}

Now we use the previous results to prove the first of our main theorems, \ref{fattree} and hence \ref{bananatree}. We need the following.
Let $(T,r)$ be a rooted tree. 
A path in $T$ of positive length joining some $u\in L(T)$ to some ancestor of $u$ is called a {\em limb} of $(T,r)$, and we call
$u$ its {\em leaf} and $v$ its {\em start}.
Let $\mathcal{T}=(T_q:q\in Q)$ be a family of limbs in $(T,r)$,
and let $k\ge 1$ be an integer. 
We make a graph $J$ with vertex set $Q$ as follows. We say that distinct $q_1,q_2\in Q$ are adjacent in $J$ 
if there are at least $k$ vertices $w$ of $T$, such that $w$ belongs to the interiors of $T_{q_1}$ and $T_{q_2}$, and $w$ is not 
a vertex of $T_q$ for any $q\in Q\setminus \{q_1,q_2\})$.
We call $J$ the {\em $k$-overlap graph} of $\mathcal{T}$. 
It is easy to see that $J$ must be a forest. More important for us
is the converse; that 
\begin{thm}\label{getforest}
For every forest $J$ and every $k\ge 1$, there is a rooted tree $(T,r)$ and a family of limbs 
$\mathcal{T}$ in $T$ such that the $k$-overlap graph of $\mathcal{T}$ is isomorphic to $J$, and no two members of $\mathcal{T}$
share an end.
\end{thm}
The proof is straightforward and we omit it.

Let us say a {\em banana} is a graph formed by the union of a nonempty set of paths each of positive length, 
all with the same ends ($s,t$ say)
and otherwise disjoint, and its {\em thickness} is the number of these paths. 
Its {\em length} is the minimum length of its constituent paths.
We call $s,t$ {\em ends} of the banana. 
By a {\em banana in $G$} we mean an induced subgraph of $G$ that is a banana.

\begin{thm}\label{getbanana}
Let $(T,r,((X_v,Y_v):v\in V(T)\setminus L(T)), (C_v:v\in V(T)))$
be a fruitful pineapple tree in $G$, and 
let $W\subseteq V(T)$.
Let
$T_1, T_2$ be limbs of $(T,r)$, with distinct leaves, and such that
for $i = 1,2$, $W$ is a subset of the interior of $T_i$.
Let $T_i$ have leaf $u_i$ and start $v_i$, and
let $x_i\in X_{v_i}^{u_i}$, where $x_1,x_2$ are nonadjacent.
There is a banana $B$ in $G$, with ends $x_1,x_2$, thickness $|W|$, and interior a subset of
$\bigcup_{w\in W} Y_w$.
\end{thm}
\Proof
Let $w\in W$. For $i = 1,2$, $(u_i,v_i,w)$ is an aligned triple, and so $x_i$ has a neighbour in $Y_w$, 
since the pineapple tree is fruitful.
Since $G[Y_w]$
is connected, there is an induced path, $P_w$ say, between $x_1, x_2$ with interior in $Y_w$. 
The union over all $w\in W$ of the 
paths $P_w$ makes the desired banana $B$.~\bbox

Two bananas $B_1,B_2$ in $G$ are {\em orthogonal} if 
every vertex in $V(B_1\cap B_2)$ is an end of both bananas, and there is at most one such vertex, and every edge of $G$
between $V(B_1)$ and $V(B_2)$ is incident with a common end of $B_1,B_2$.  A {\em banana tree} is a graph obtained from a tree $T$
by replacing each edge $uv$ by a banana with ends $u,v$, such that these bananas are orthogonal. To prove
\ref{fattree}, we need to prove that every multigraph obtained by fattening the edges of a tree is widespread; and part of
``widespread'' involves proving that for every subdivision $J$ of such a multigraph, there is a subdivision of $J$ which is
present as an induced subgraph. But such a graph $J$ is just a banana tree, so \ref{fattree} is implied by the following:

\begin{thm}\label{fattree3}
For every banana tree $J$ and 
for all $\tau\ge 0$ there exists $c\ge 0$ such that every $J$-subdivision-free graph $G$ with
$\chi^{2|V(J)|+7}(G)\le \tau$ satisfies $\chi(G)\le c$.
\end{thm}
\Proof Let $n=|V(J)|$. We may assume that $n\ge 3$, since otherwise the result is trivial. 
Since $J$ is a banana tree, it is obtained from some tree $S$ by substituting bananas for its edges.
By \ref{getforest}, we may choose a rooted tree $(T,r)$ such that there is a family $\mathcal{T}$ 
of limbs in $(T,r)$ with $n$-overlap graph isomorphic to $S$, and no two members of $\mathcal{T}$ share an end.
Choose $\nu>0$ such that there is an induced subgraph
of $K_{\nu}^1$ which is a subdivision of $J$; and let
$c'$ satisfy \ref{summary} with $c,d$ replaced by $0,n$.
We claim that setting $c=c'$ satisfies the conclusion of the theorem. For let $G$ be a $J$-subdivision-free graph $G$ with
$\chi^{2n+7}(G)\le \tau$. We must show that $\chi(G)\le c'$. Suppose not. Now
$G$ is $K_{\nu}^1$-subdivision-free, and so by \ref{summary}
(setting $Z=V(G)$) there is a fruitful pineapple tree
$$(T,r,((X_v,Y_v):v\in V(T)\setminus L(T)), (C_v:v\in V(T)))$$
in $G$, such that
for all incomparable $u,v\in L(T)$, the $G$-distance between $C_u,C_v$
is at least $n+1$.
For each $u\in L(T)$, $C_u$ is nonempty from the definition of a pineapple tree.
Let $\mathcal{T} =(T_q:q\in Q)$, and for each $q\in Q$ let $T_q$ have leaf $u_q$ and start $v_q$; and choose
$x_q\in X_{v_q}^{u_q}$ (this is possible since $X_{v_q}^{u_q}$ covers $C_{u_q}$). 
There is an isomorphism from the $n$-overlap graph of $\mathcal{T}$ to $S$; let the corresponding
bijection from $Q$ onto $V(S)$ map $q$ to $s_q$ for each $q\in Q$.
Now let $s_qs_{q'}$
be an edge of $S$. Then from the definition of the $n$-overlap graph, there are at least $n$ vertices of $T$ that belong to the 
interiors of the limbs $T_q, T_{q'}$ and do not belong to any of the other members of $\mathcal{T}$.
Let $W_e$ be a set of $n$ such vertices. 
By \ref{getbanana}, there is a banana $B_e$ with ends $x_q,x_{q'}$, with thickness $n$, 
such that its interior is a subset of $\bigcup_{w\in W_e}Y_w$. The $G$-distance between $x_q,x_{q'}$ is at least $n-1$, because
the $G$-distance between $C_{u_q}$ and $C_{u_{q'}}$ is at least $n+1$; so $B_e$ has length at least $n-1$.

Since the vertices in $W_e$ do not belong to any other member of $\mathcal{T}$, it follows that for all distinct edges
$e,f$ of $S$, the 
bananas $B_e,B_f$ are disjoint except for their ends. Since each banana has thickness $n$ and length at least $n-1$, 
it follows that
the vertices $x_1\l x_k$ and all the bananas $B_e$ make a subgraph $H$ of $G$ which has an induced subgraph that is a 
subdivision of $J$. We claim that $H$ is itself an induced subgraph of $G$. For suppose not, and
let $a,b\in V(H)$ be distinct and adjacent in $G$ and not adjacent in $H$. Since the vertices $x_q\;(q\in Q)$ pairwise have $G$-distance
at least $n-1$ and hence are nonadjacent (since $n\ge 3$), we may assume that $a$ belongs to the interior of some banana $B_e$
say, and hence $a\in Y_w$ for some $w$. 

From the definition of a pineapple tree, the only vertices of the pineapple tree with neighbours in $Y_w$
belong to $Y_w$, to $X_w$, or to $X_v$ for some ancestor $v$ of $w$; and if a vertex in $X_v^u$ has a neighbour in $Y_w$
where $v$ is an ancestor of $w$ and $u\in L(T)$, then $u$ is a descendant of $w$. Consequently the only vertices of $H$
with neighbours in $Y_w$ belong to $Y_w\cup \{x_q,x_{q'}\}$ where $e=s_qs_{q'}$; and so $b$ belongs to this set, and hence to 
$P_e^w$. Since $P_e^w$ is induced, it follows that $a,b$ are adjacent in $H$, a contradiction. This proves that
$H$ is induced. Consequently there is an induced subgraph of $G$ isomorphic to a subdivision of $J$, a contradiction.
This proves \ref{fattree3}.~\bbox

\section{Fattening a cycle}

In this section we prove \ref{cycle}. In the proof of \ref{fattree3} we made use of the overlap graph, which exploited  vertices
of the tree that only belonged to two of the selected limbs. For \ref{cycle} we will again apply \ref{gettree}, 
but now we need to use vertices of the tree that belong to more than two limbs.
This will still give
us bananas, but we have less control over which pairs the bananas join. We have the following variant of \ref{getbanana}.

\begin{thm}\label{bigshare}
Let $(T,r,((X_v,Y_v):v\in V(T)\setminus L(T)), (C_v:v\in V(T)))$
be a fruitful pineapple tree in $G$, such that 
for all incomparable $u,v\in L(T)$, the $G$-distance between $C_u,C_v$ 
is at least $5$.
Let $k,n\ge 1$ be integers, and let $W\subseteq V(T)$ with $|W|>(n-1)k(k-1)/2$.
Let 
$(T_q:q\in Q)$ be a family of limbs of $(T,r)$, with $|Q|=k$, such that their leaves are all different, and such that 
$W$ is a subset of the interior of $T_q$ for each $q\in Q$.
For each $q\in Q$ let $T_q$ have leaf $u_q$ and start $v_q$, and 
let $x_q\in X_{v_q}^{u_q}$. 
For every partition of $Q$ into two nonempty sets
$I,J$, there exist $q\in I$, $q'\in J$, and a banana $B$ in $G$, with ends $x_q,x_{q'}$, thickness $n$ and interior a subset of
$\bigcup_{w\in W} Y_w$, such that there is no edge between $V(B)$ and $\{x_{q''}:q''\in Q\setminus \{q,q'\}\}$.
\end{thm} 
\Proof 
Since $|Q|=k$, we may assume that $Q=\{1\ll k\}$.
We see first that for all distinct $q,q'\in Q$, the $G$-distance between $x_q,x_{q'}$ is at least three, because they have
neighbours in $C_{u_q}, C_{u_{q'}}$ respectively, and the $G$-distance between $C_{u_q}, C_{u_{q'}}$ is at least five, since
$u_q\ne u_{q'}$. It follows
that the vertices $x_1\ll x_k$ are pairwise nonadjacent and no two have a common neighbour. Let $I,J$ be complementary nonempty subsets 
of $Q$, and let $w\in W$. For each $i\in I$ and $j\in J$, since $x_i,x_j$ both have neighbours in $Y_w$ and $G[Y_w]$
is connected, there is a path between $x_i, x_j$ with interior in $Y_w$. Choose $i,j$ and the path such that this path 
($P_w$ say) is as short
as possible. It follows that no other member of $\{x_1\ll x_k\}$ has a neighbour in $P_w$, since 
 $x_1\ll x_k$ are pairwise nonadjacent and no two have a common neighbour. Let $f_w=(i,j)$ where $i<j$. Since there are only
$k(k-1)/2$ possibilities for the pair $(i,j)$ (in fact fewer, since we insists that $i\in I$ and $j\in J$), there exist
at least $n$ values of $w$ where the pairs $f_w$ are all the same, equal to $(q,q')$ say. But then the union of the corresponding
paths $P_w$ makes the desired banana $B$.~\bbox

The goal of the section is to prove \ref{cycle}, which is implied by the following:
\begin{thm}\label{cycle3}
Let $J$ be obtained from a cycle of length $m$ by substituting bananas for all except one of its edges.
Then for all $\tau\ge 0$ there exists $c\ge 0$ such that  
if $G$ is a $J$-subdivision-free graph
with $\chi^{2|V(J)|+7}(G)\le \tau$, then $\chi(G)\le c$.
\end{thm}
\Proof
Let $n=\max(|V(J)|,5)$.
Let $(S,r)$ be a uniform $2$-ary tree of height $m$. For each vertex $z\in V(S)$, we say its {\em height} is the $S$-distance
from $z$ to a vertex in $L(S)$.
Let $(T,r)$ be obtained from $(S,r)$ by replacing each edge $e\in E(S)$ 
by a path
$P_e$ of length $2n^3$. (Thus $V(S)\subseteq V(T)$.) 
Choose $\nu$ such that $K^1_{\nu}$ contains $J$; and 
choose $c'$ such that \ref{summary} holds, taking $d=n+1$ and $c=0$. We claim that setting $c=c'$ satisfies the conclusion of the theorem.

Let $G$ be a $J$-subdivision-free graph with 
$\chi^{2|V(J)|+7}(G)\le \tau$, and suppose that $\chi(G)> c$. 
Let $h\ge 1$, and let $x_1\ll x_h\in V(G)$ be distinct. 
A {\em banana path on $(x_1\l x_h)$} is a sequence $(B_1\ll B_{h-1})$ of 
pairwise orthogonal bananas in $G$, each of thickness $n-1$,
such that $B_i$ has ends
$x_i, x_{i+1}$ for $1\le i\le h-1$. (If $h=1$, the null sequence of bananas counts as a banana path.)
Its {\em interior} is the union of the interiors of $B_1\l B_{h-1}$. (Thus none of $x_1\ll x_h$
belongs to the interior.)

Since $G$ is $K_{\nu}^1$-subdivision-free, by \ref{summary} there is a fruitful pineapple tree
$$(T,r,((X_v,Y_v):v\in V(T)\setminus L(T)), (C_v:v\in V(T)))$$
in $G$, such that
for all incomparable $u,v\in L(T)$, the $G$-distance between $C_u,C_v$
is at least $n+2$.
For each $u\in L(T)$ choose $x^u\in X_r^u$. 
Since the set $C_u\;(u\in L(T))$
pairwise have $G$-distance at least $n+2$, it follows that the vertices 
$x^u\;(u\in L(T))$ pairwise have $G$-distance at least $n$.
For each edge $e=z_1z_2 \in E(S)$, where $z_1$ is the parent of $z_2$ in $S$, let $T_e$ be the union of 
$V(P_e)\setminus \{z_1\}$ and the set of       
descendants of $z_2$ in $T$. Let $L_e=T_e\cap L(T)$, and let $M_e=\{x^u:u\in L_e)\}$. Thus $M_e\subseteq X_r$ for 
every $e\in E(S)$.
\\
\\
(1) {\em Let $e=z_1z_2 \in E(S)$, where $z_1$ is the parent of $z_2$ in $S$, and let $h$ be the height of $z_1$. Then 
there exist distinct vertices $x_1\l x_h\in M_{e}$, such that there is a banana path on $(x_1\l x_h)$ with interior
a subset of $\bigcup_{v\in T_e\setminus L_e}Y_v$.}
\\
\\
To prove this we proceed by induction on $h$. If $h=1$, then the result is true since then $z_2\in L(T)$ and $x^{z_2}$
satisfies the requirement. Thus we may assume that $h>1$ and the result holds for $h-1$.
Let $z_3,z_4$ be the children of $z_2$ in $S$, and let $f,g$ be the edges $z_2z_3$ and $z_2z_4$ respectively.
From the inductive hypothesis, there exist distinct $a_1\l a_{h-1}\in M_{z_3}$ such that there is a banana path
on $(a_1\l a_{h-1})$ with interior in $\bigcup_{v\in T_f\setminus L_f}Y_v$, and there exist $b_1\l b_{h-1}\in M_{z_4}$
similarly. It follows that there are no edges between any banana of the first banana path and any banana of the second, 
from the definition of a pineapple tree. Moreover, each of $a_1\l a_{h-1}$ has $G$-distance at least $n$ from each of
$b_1\l b_{h-1}$. Let $W_e$ denote the set of vertices of $T$ in the interior of $P_e$.
Now each $a_i$ is adjacent to a member of $\bigcup_{u\in L(T)} C_u$, and the interior of the corresponding limb includes $W_e$,
and the same for each $b_j$. Since $|W_e|\ge 2n^3-1>(n-1)2h(2h-1)/2$,
\ref{bigshare} implies that there exist $a_i\in \{a_1\l a_{h-1}\}$ and $b_j\in \{b_1\l b_{h-1}\}$, and a banana $B$ with ends
$a_i,b_j$ and thickness $n$ and interior a subset of
$\bigcup_{w\in W_e} Y_w$, such that no other member of $\{a_1\l a_{h-1},b_1\l b_{h-1}\}$ has a neighbour in $V(B)$.
By reversing the sequence $(a_1\l a_{h-1})$ if necessary, we may assume that $i\ge h/2$, and similarly $j\le h/2$.
Now no vertex of the interior of $B$ belongs to or has a neighbour in $\bigcup_{v\in T_f\setminus L_f}Y_v$, and the same for
$\bigcup_{v\in T_g\setminus L_g}Y_v$; and so there is a banana path on
$(a_1\l a_i,b_j\l b_{h-1})$ with interior in 
$$\bigcup_{v\in T_f\setminus L_f}Y_v\cup \bigcup_{v\in T_g\setminus L_g}Y_v\cup \bigcup_{v\in W_e}Y_v\subseteq \bigcup_{v\in T_e\setminus L_e}Y_v.$$
This proves (1).

\bigskip

In particular, since $r$ has height $m$, from (1) applied to some edge $e=rs$ of $S$ incident with $r$, 
there exist $m$ vertices $x_1\l x_m\in M_{e}$ such that there is a banana path on $(x_1\l x_m)$ with interior
a subset of $\bigcup_{v\in T_e\setminus L_e}Y_v$. Now choose a path between $x_1,x_m$ with interior in $Y_r$ such that
$x_2\l x_{m-1}$ have no neighbours in it (this is possible since the pineapple $(X_r,Y_r)$ is levelled). Adding this to the 
banana path gives a subdivision of $J$. This proves \ref{cycle3}.~\bbox

\section{The fat triangle}

Now we turn to the third of our theorems, \ref{triangle}. We need a lemma, as follows.

\begin{thm}\label{internalroute}
Let $\rho\ge 4$, $\tau\ge 0$, and $n\ge 0$,
let $G$ be a graph with $\chi^{\rho}(G)\le \tau$, and let $X,Z\subseteq V(G)$ be disjoint, such that $X$ covers $Z$ 
and $\chi(Z)>n\tau$. 
Then there exist $x_1\l x_n\in X$ with the following properties:
\begin{itemize}
\item $x_1\l x_n$ pairwise have $G$-distance at least $\rho$; and
\item for all distinct $i,j\in \{1\l n\}$, there is a path between $x_i,x_j$ with interior in $Z$, such that
no other vertex in $\{x_1\l x_n\}$ has a neighbour in this path.
\end{itemize}
\end{thm}
\Proof
By a {\em $\{1\l k\}$-colouring} of a graph we mean a proper colouring using $\{1\l k\}$ as the set of colours.
If $X_1,X_2\subseteq V(G)$ with $X_1\cap X_2 = \emptyset$, and $\kappa_i$ is a colouring of $G[X_i]$ for $i = 1,2$, we say they
are {\em compatible} if their union in the natural sense is a colouring of $G[X_1\cup X_2]$.
For each $x\in X$, the subgraph $G_x$ induced on the set of neighbours of $x$ in $Z$
has chromatic number at most $\tau$; choose some $\{1\l \tau\}$-colouring $\kappa_x$ of $G_x$ for each such $v$.
Now choose $k\in \{0\ll n\}$, $x_1\l x_k\in X$ and  $C\subseteq Z$, with the following properties:
\begin{itemize}
\item $G[C]$ is connected, and $x_1\l x_k$ have no neighbours in $C$;
\item $x_1\l x_k$ pairwise have $G$-distance at least $\rho$;
\item no $\{1\l n\tau\}$-colouring of $G[C]$ is compatible with each of the colourings $\kappa_{x_i}\; (1\le i\le k)$; and
\item subject to these conditions, $C$ is minimal.
\end{itemize}
(This is possible, since taking $k=0$ and taking $C$ to be the vertex set of a component of $G[Z]$ with maximum chromatic number satisfies all bullets
except the last.) If some $G_{x_i}$ contains no vertex with a neighbour in $C$ then we may remove $x_i$ from the list $x_1\l x_k$;
so we may assume that $x_1\l x_k$ each have a neighbour in $Z$ which has a neighbour in $C$.
Now $k\le n$, and if $k = n$ then the theorem holds, so we may suppose for a contradiction that $k<n$.
Since only colours $1\l \tau$ are used by the colourings $\kappa_1\l \kappa_k$,
it follows that
$\chi(C)>(n-1)\tau\ge k\tau$; and so there exists $v\in C$ with $G$-distance more than $\rho$ from each of $x_1\l x_k$. Choose
$x_{k+1}\in X$ adjacent to $v$. Thus $x_{k+1}$ has $G$-distance at least $\rho$ from each of $x_1\l x_k$.
Let $C'$ be the set of vertices in $C$ nonadjacent to $x_{k+1}$.
Since $\kappa_{x_{k+1}}$
is compatible with each of the colourings $\kappa_{x_i}\; (1\le i\le k)$ (because $x_{k+1}$ has $G$-distance at least four from each
of $x_1\l x_k$), it follows that
no $\{1\l n\tau\}$-colouring of $G[C']$ is compatible with each of the colourings $\kappa_{x_i}\; (1\le i\le k+1)$.
Consequently there is a component $C''$ of $G[C']$ such that no $\{1\l n\tau\}$-colouring of $G[C'']$
is compatible with each of the colourings $\kappa_{x_i}\; (1\le i\le k+1)$.
But this contradicts the minimality of $C$. Hence $k = n$. This proves \ref{internalroute}.~\bbox

We deduce the following, which implies \ref{triangle}:

\begin{thm}\label{triangle2}
Let $H$ be the multigraph obtained from $K_3$ by fattening two of its edges and replacing the third by two parallel edges, 
and let $J$ be a subdivision of $H$. Let $n=|V(J)|$.
For all $\tau\ge 0$
there exists $c\ge 0$ such that every $J$-subdivision-free graph $G$ with
$\chi^{2n+7}(G)\le \tau$ satisfies $\chi(G)\le c$.
\end{thm}
\Proof
It follows that $n\ge 3$.
Let $(T,r)$ be the rooted tree where $T$ is a path of length $3n$ and $r$ is one end of $T$. Let the vertices of $T$ in order be
$t_0\l t_{3n}$ where $r=t_0$; thus $L(T) = \{t_{3n}\}$. We write $u=t_{3n}$.
Choose $\nu$ such that $K_{\nu}^1$ contains $J$.
Choose $c'$ to satisfy \ref{summary} with $c=3\tau$ and $d=n$; we claim that setting $c=c'$ satisfies the conclusion of the theorem.
For let $G$ be $J$-subdivision-free and hence $K_{\nu}^1$-subdivision-free, with $\chi^{2n+7}(G)\le \tau$, and suppose that
$\chi(G)>c'$, By \ref{summary},
there is a fruitful pineapple tree
$$(T,r,((X_v,Y_v):v\in V(T)\setminus L(T)), (C_v:v\in V(T)))$$
in $G$, such that
$\chi(C_u)>3\tau$.
For $0\le i<3n$ let $X_i=X_{t_i}$ and $Y_i=Y_{t_i}$. We may assume that every vertex in $X_i$ has a neighbour in $C_u$,
because any other vertices in $X_i$ may be removed from $X_i$ (thus $X_{t_i} = X_{t_i}^{t_{3n}}$ in the earlier notation).
Consequently for all $i,j$ with $0\le i<j<3n$, every vertex in $X_i$ has a neighbour in $Y_j$.

By \ref{internalroute}, there exist
$x,x',x''\in X_0$, pairwise at $G$-distance at least $n$,
such that every two of them are joined by a path with interior in $C_u$ in which the third has no neighbours.
Let $R(xx')$ be the set of $i\in \{1\l 3n-1\}$ such that there is a path
between $x,x'$ with interior in $Y_i$ containing no neighbour of $x''$; and define $R(xx'')$ and $R(x'x'')$ similarly.
Since 
$x,x',x''$ pairwise have $G$-distance at least $n\ge 3$, and $G[Y_i]$ is connected and $x,x',x''$ 
have neighbours in $Y_i$, 
each value of $i\in \{1\l 3n-1\}$ belongs to at least two of $R(xx'), R(xx''), R(x'x'')$. Consequently there exists 
$I\subseteq \{1\l 3n-1\}$ with $|I|=n$ such that $I$ is a subset of one of $R(xx'), R(xx''), R(x'x'')$, say $R(xx'')$; and since
$|\{1\l 3n-1\}\setminus I|=2n-1$, there exists $J\subseteq \{1\l 3n-1\}\setminus I$ with $|J|=n$ such that $J$ is a subset of one
of $R(xx'), R(x'x'')$, say $R(x'x'')$. By \ref{getbanana}, there is a banana of thickness $n$, with ends
$x,x''$, and with interior a subset of $\bigcup_{i\in I}Y_i$; and it has length at least $n$, since the $G$-distance between $x,x''$
is at least $n$. Similarly there is a banana with ends $x',x''$
with interior a subset of $\bigcup_{j\in J} Y_j$. These two bananas are orthogonal. 

To obtain a subdivision
of $J$, we need to add to this union two paths joining $x,x'$; and we will
obtain these, one with interior in $C_u$ via \ref{internalroute},
and one with interior in $Y_0$.
The first is immediate from the definition of $x,x',x''$. For the second we use the fact that $(X_0,Y_0)$ is a levelled pineapple. 
Let $z_0\in Y_0$ such that for some $k$, every vertex in $Y_0$ is joined to $z_0$ by a path of $G[Y_0]$ of length less than $k$,
and there is no path in $G[X_0\cup Y_0]$ from $z_0$ to $X_0$ of length less than $k$. For $0\le i\le k$, let $L_i$
be the set of vertices in $X_0\cup Y_0$ with $G[X_0\cup Y_0]$-distance $i$ from $z_0$.
Thus $Y_0 = L_0\cup\cdots\cup L_{k-1}$ and $X_0=Y_k$. Now $x,x'$ both have neighbours in $L_{k-1}$, say $y,y'$ respectively, and
since the $G$-distance  between $x,x'$ is at least $n$, it follows that $2k\ge n$, and in particular $k>1$.
Since $G[L_0\cup\cdots\cup L_{k-2}]$ is connected and $y,y'$ both have neighbours in it, there is an induced path between $y,y'$
with interior in $L_0\cup\cdots\cup L_{k-2}$, and which consequently contains no neighbours of $x''$. 
Adding the edges $xy$ and $x'y'$ to this path gives the required path from $x$ to $x'$. 
The subgraph consisting of the two bananas and these two paths is induced, and isomorphic to a subdivision of $H$.
Since $x,x',x''$ pairwise have distance at least $n$, all these paths between
pairs of $x,x',x''$ have length at least $n$; and so this same subgraph is also isomorphic to a subdivision of $J$, 
a contradiction. This proves \ref{triangle2}.~\bbox

\section{Bigger widespread graphs}

The same methods can be combined to prove that more complicated graphs are widespread. For instance, in the proof of \ref{cycle}, all 
the limbs we used started from the root $r$, and all the limbs eventually become disjoint. We are free to make
the tree $T$ bigger by adding more vertices to its leaves, and extend the old limbs further to make new limbs, 
and add more limbs meeting the old paths just in their new sections.
By this process we can make not just one cycle as in \ref{cycle}, but any multigraph each of whose blocks is such a cycle. 
We omit the details.

Can we make more 2-connected widespread graphs? Here is one construction. Take a path with vertices $v_1\c v_k$ in order 
where $k\ge 4$, fatten each edge, and add two more vertices $a,b$ and edges $av_1,av_2, bv_{k-1},bv_k$ and $ab$, making
a multigraph $H$. We claim:

\begin{thm}\label{oddgraph}
$H$ is widespread.
\end{thm}
\Proof
We merely sketch the proof, since the result is such an oddity. Let $J$ be a subdivision of $H$.
Let us proceed as in the proof of \ref{cycle3}, with a subdivided
$2$-ary tree $T$; but apply \ref{summary} to this tree with $c$ larger than zero, 
large enough that \ref{internalroute} can be applied. For each $u\in L(T)$, choose three vertices $x_u,y_u,z_u\in X_r^u$,
such that for every two of them there is a path between them with interior in $C_u$ containing no neighbour of the third; and
choose $x_u,y_u,z_u$ with $G$-distance at least $n+2$. Now because the limb of $T$ from $r$ to $u$
has a final section $W$ consisting of many vertices $w$ each with only one child, these vertices are incomparable with the other
leaves of $T$, and so there are two orthogonal bananas $B_u,B'_u$ in $G$, both with interior in 
the union of $Y_w\; (w\in W)$, and each with both ends in $\{x_u,y_u,z_u\}$ (and joining distinct pairs from this set). 
This defines a banana path of length two.
Now we apply the method of \ref{cycle3}; we generate longer and longer banana paths, starting from the ones we just made of length two.
The procedure of \ref{cycle3} has the convenient feature that the first and last banana of every banana path it generates
is a banana of one of the initial banana paths of length two. So we may assume that we generate a $(k-1)$-term banana path
where the first banana is $B_u$ and the last is $B_{u'}$ for some $u,u'\in L(T)$. Let $B_u$ have ends $x_u,y_u$ say.
By choosing one path from the banana $B_u'$ (not to be confused with $B_{u'}$), joining $z_u$ with one of $x_u,y_u$,
and choosing one path via \ref{internalroute} joining $z_u$ with the other of $x_u,y_u$, we obtain an induced path from $x_u$
to $y_u$ in which $z_u$ is an internal vertex. Now do the same thing for $u'$, and then add a path joining $z_u,z_{u'}$
with interior in $Y_r$. This provides the induced subgraph which is a subdivision of $J$.~\bbox

\end{document}